\documentclass[11pt]{article}
\usepackage[]{amsmath,amssymb}

\newtheorem{theorem}{Theorem}[section]
\newtheorem{lemma}[theorem]{Lemma}
\newtheorem{proposition}[theorem]{Proposition}
\newtheorem{corollary}[theorem]{Corollary}
\newtheorem{exAux}[theorem]{Example}

\newtheorem{Def}[theorem]{Definition}
\newenvironment{definition}{\begin{Def} \rm}{\end{Def}}
\newtheorem{Note}[theorem]{Note}
\newenvironment{note}{\begin{Note} \rm}{\end{Note}}
\newtheorem{Rem}[theorem]{Remark}

\newtheorem{Ass}[theorem]{Assumption}

\newenvironment{proof}{\medskip\noindent{\bf Proof.\ }}{\qed\medskip}
\newenvironment{proofof}[1]{\medskip\noindent{\bf Proof  of {#1}.\ 
}}{\qed\medskip}
\newcommand{\qed}{\hfill\mbox{$\Box$\qquad\qquad}}

\newcommand{\Mat}[1]{\text{\rm Mat}_{#1}(\mathbb{K})}

\newcommand{\T}{\Upsilon}

\begin{document}
\thispagestyle{empty}

\begin{center}
\LARGE \bf
\noindent
Linear transformations that are tridiagonal with respect to
both eigenbases of a Leonard pair
\end{center}

\smallskip

\begin{center}
\Large
Kazumasa Nomura and Paul Terwilliger
\end{center}

\smallskip

\begin{quote}
\small 
\begin{center}
\bf Abstract
\end{center}
Let $\mathbb{K}$ denote a field, and let $V$ denote a vector
space over $\mathbb{K}$ with finite positive dimension.
We consider a pair of linear transformations 
$A : V \to V$ and $A^* : V \to V$ that satisfy (i) and (ii) below:
\begin{itemize}
\item[(i)] There exists a basis for $V$ with respect to which the
matrix representing $A$ is irreducible tridiagonal and the matrix
representing $A^*$ is diagonal.
\item[(ii)] There exists a basis for $V$ with respect to which the
matrix representing $A^*$ is irreducible tridiagonal and the matrix
representing $A$ is diagonal.
\end{itemize}
We call such a pair a {\em Leonard pair} on $V$.
Let $\cal X$ denote the set of linear transformations $X:V \to V$
such that the matrix representing $X$ with respect to 
the basis (i) is tridiagonal and 
the matrix representing $X$ with respect to the basis (ii) is tridiagonal.
We show that $\cal X$ is spanned by 
\[
    \text{$I$, $A$, $A^*$, $AA^*$, $A^*A$},
\]
and these elements form a basis for $\cal X$
provided the dimension of $V$ is at least $3$.
\end{quote}

\section{Leonard pairs}

We begin by recalling the notion of a Leonard pair.
We will use the following terms.
A square matrix $X$ is said to be {\em tridiagonal}
whenever each nonzero entry lies on either the diagonal, the subdiagonal,
or the superdiagonal. Assume $X$ is tridiagonal.
Then $X$ is said to be {\em irreducible}
whenever each entry on the subdiagonal is nonzero and each entry on
the superdiagonal is nonzero.
We now define a Leonard pair.
For the rest of this paper $\mathbb{K}$ will denote a field.

\medskip

\begin{definition}  \cite{T:Leonard}  \label{def:LP}
Let $V$ denote a vector space over $\mathbb{K}$ with finite positive
dimension.
By a {\em Leonard pair} on $V$ we mean an ordered pair $(A,A^*)$,
where $A:V \to V$ and $A^*:V \to V$ are linear transformations
that satisfy (i) and (ii) below:
\begin{itemize}
\item[(i)] There exists a basis for $V$ with respect to which the
matrix representing $A$ is irreducible tridiagonal and the matrix
representing $A^*$ is diagonal.
\item[(ii)] There exists a basis for $V$ with respect to which the
matrix representing $A^*$ is irreducible tridiagonal and the matrix
representing $A$ is diagonal.
\end{itemize}
\end{definition}

\medskip

\begin{note}
It is a common notational convention to use $A^*$ to represent the
conjugate-transpose of $A$. We are {\em not} using this convention.
In a Leonard pair $(A,A^*)$ the linear transformations $A$ and
$A^*$ are arbitrary subject to (i) and (ii) above.
\end{note}

\medskip

We refer the reader to 
\cite{H},
\cite{N:aw},
\cite{NT:balanced}, \cite{NT:formula}, \cite{NT:det}, \cite{NT:mu},
\cite{P}, \cite{T:sub1}, \cite{T:sub3}, \cite{T:Leonard},
\cite{T:24points}, \cite{T:conform}, \cite{T:intro},
\cite{T:intro2}, \cite{T:split}, \cite{T:array}, \cite{T:qRacah},
\cite{T:survey}, \cite{TV}, \cite{V}
for background on Leonard pairs.
We especially recommend the survey \cite{T:survey}.
See \cite{AC}, \cite{AC2}, \cite{ITT}, \cite{IT:shape},
\cite{IT:uqsl2hat}, \cite{IT:non-nilpotent}, \cite{ITW:equitable}, 
\cite{N:refine}, \cite{N:height1},
\cite{T:qSerre}, \cite{T:Kac-Moody} for related topics.

\section{Leonard systems}

When working with a Leonard pair, it is convenient to consider a closely
related object called a {\em Leonard system}. 
To prepare for our definition
of a Leonard system, we recall a few concepts from linear algebra.
Let $d$ denote a nonnegative integer and let
$\Mat{d+1}$ denote the $\mathbb{K}$-algebra consisting of all $d+1$ by
$d+1$ matrices that have entries in $\mathbb{K}$. 
We index the rows and  columns by $0, 1, \ldots, d$. 
For the rest of this paper, let
$\cal A$ denote a $\mathbb{K}$-algebra isomorphic to $\Mat{d+1}$,
and let $V$ denote a simple $\cal A$-module. We remark that $V$ is unique
up to isomorphism of $\cal A$-modules, and that $V$ has dimension $d+1$.
Let $v_0, v_1, \ldots, v_d$ denote a basis for $V$.
For $X \in {\cal A}$ and $Y \in \Mat{d+1}$, we say 
{\em $Y$ represents $X$ with respect to} $v_0, v_1, \ldots, v_d$
whenever $X v_j = \sum_{i=0}^d Y_{ij}v_i$ for $0 \leq j \leq d$.
For $A \in \cal A$ we say $A$ is {\em multiplicity-free}
whenever it has $d+1$ mutually distinct eigenvalues in $\mathbb{K}$. 
Assume $A$ is multiplicity-free. 
Let $\theta_0, \theta_1, \ldots, \theta_d$ denote an ordering 
of the eigenvalues of $A$, and for $0 \leq i \leq d$ put
\begin{equation}        \label{eq:defEi}
    E_i = \prod_{\stackrel{0 \leq j \leq d}{j\neq i}}
             \frac{A-\theta_j I}{\theta_i - \theta_j},
\end{equation}
where $I$ denotes the identity of $\cal A$. 
We observe
(i) $AE_i = \theta_i E_i$ $(0 \leq i \leq d)$;
(ii) $E_i E_j = \delta_{i,j} E_i$ $(0 \leq i,j \leq d)$;
(iii) $\sum_{i=0}^{d} E_i = I$;
(iv) $A = \sum_{i=0}^{d} \theta_i E_i$.
Let $\cal D$ denote the subalgebra of $\cal A$ generated by $A$.
Using (i)--(iv) we find the sequence $E_0,E_1,\ldots,E_d$
is a basis for the $\mathbb{K}$-vector space $\cal D$.
We call $E_i$ the {\em primitive idempotent} of $A$ associated with
$\theta_i$. 
It is helpful to think of these primitive idempotents as follows.
Observe 
\begin{equation}           \label{eq:decomp}
 V=E_0V+E_1V+\cdots+E_dV  \qquad \text{(direct sum)}.
\end{equation}
For $0 \leq i \leq d$, $E_iV$ is the (one dimensional) eigenspace of $A$
in $V$ associated with the eigenvalue $\theta_i$, and $E_i$ acts on $V$
as the projection onto this eigenspace.
We remark that the $\mathbb{K}$-vector space $\cal D$ has basis
$\{A^i\,|\, 0 \leq i \leq d\}$ and satisfies
${\cal D}=\{X \in {\cal A}\,|\, AX=XA\}$.

\medskip

By a {\em Leonard pair in $\cal A$} we mean an ordered pair of elements
taken from $\cal A$ that act on $V$ as a Leonard pair in the sense of
Definition \ref{def:LP}.
We now define a Leonard system.

\medskip

\begin{definition}  \cite{T:Leonard}     \label{def:LS}   \samepage
By a {\em Leonard system} in $\cal A$ we mean a sequence
\[
   (A; \{E_i\}_{i=0}^d; A^*; \{E^*_i\}_{i=0}^d)
\]
that satisfies (i)--(v) below.
\begin{itemize}
\item[(i)] Each of $A$, $A^*$ is a multiplicity-free element in $\cal A$.
\item[(ii)] $E_0, E_1, \ldots, E_d$ is an ordering of the
   primitive idempotents of $A$.
\item[(iii)] $E^*_0, E^*_1, \ldots, E^*_d$ is an ordering of the
   primitive idempotents of $A^*$.
\item[(iv)] For $0 \leq i,j \leq d$, 
\begin{equation}                    \label{eq:Astrid}
   E_i A^* E_j =
    \begin{cases}  
        0 & \text{\rm if $|i-j|>1$},  \\
        \neq 0 & \text{\rm if $|i-j|=1$}.
    \end{cases}
\end{equation}
\item[(v)] For $0 \leq i,j \leq d$, 
\begin{equation}             \label{eq:Atrid}
   E^*_i A E^*_j =
    \begin{cases}  
        0 & \text{\rm if $|i-j|>1$},  \\
        \neq 0 & \text{\rm if $|i-j|=1$}.
    \end{cases}
\end{equation}
\end{itemize}
\end{definition}

\medskip

Leonard systems are related to Leonard pairs as follows.
Let $(A; \{E_i\}_{i=0}^d;A^*; \{E^*_i\}_{i=0}^d)$ denote a Leonard system
in $\cal A$. Then $(A,A^*)$ is a Leonard pair in $\cal A$
\cite[Section 3]{T:qRacah}.
Conversely, suppose $(A,A^*)$ is a Leonard pair in $\cal A$.
Then each of $A,A^*$ is multiplicity-free \cite[Lemma 1.3]{T:Leonard}.
Moreover there exists an ordering $E_0,E_1,\ldots,E_d$ of the
primitive idempotents of $A$, and 
there exists an ordering $E^*_0,E^*_1,\ldots,E^*_d$ of the
primitive idempotents of $A^*$, such that
$(A; \{E_i\}_{i=0}^d; A^*; \{E^*_i\}_{i=0}^d)$
is a Leonard system in $\cal A$ \cite[Lemma 3.3]{T:qRacah}.

\section{The space $\cal X$}

In this paper we consider a subspace of $\cal A$ defined as follows.

\medskip

\begin{definition}           \label{def:X}        \samepage
Let $(A; \{E_i\}_{i=0}^d; A^*; \{E^*_i\}_{i=0}^d)$
denote a Leonard system in $\cal A$.
Let $\cal X$ denote the $\mathbb{K}$-subspace of $\cal A$ consisting
of the $X \in {\cal A}$ such that both
\begin{eqnarray}
   E_i X E_j = 0  &\qquad& \text{if $\;|i-j|>1$},   \label{eq:trid} \\
   E^*_i X E^*_j = 0  &\qquad& \text{if $\;|i-j|>1$}  \label{eq:trids}
\end{eqnarray}
for $0 \leq i,j \leq d$.
\end{definition}

\medskip

We now state our main result.

\medskip

\begin{theorem}          \label{thm:main}          \samepage
Let $(A; \{E_i\}_{i=0}^d; A^*; \{E^*_i\}_{i=0}^d)$
denote a Leonard system in $\cal A$.
Then the space $\cal X$ from Definition \ref{def:X} is spanned by
\begin{equation}      \label{eq:basisX}
   I,\;A,\:A^*,\;AA^*,\;A^*A.
\end{equation}
Moreover $(\ref{eq:basisX})$ is a basis for $\cal X$ provided $d \geq 2$.
\end{theorem}

\medskip

The proof of Theorem \ref{thm:main} will be given in Section 5.

\section{The antiautomorphism $\dagger$}

Associated with a given Leonard system in $\cal A$, there is certain
antiautomorphism of $\cal A$ denoted by $\dagger$ and defined below.
Recall an {\em antiautomorphism} of $\cal A$ is an isomorphism
of $\mathbb{K}$-vector spaces $\sigma : {\cal A} \to {\cal A}$ 
such that  $(XY)^\sigma = Y^\sigma X^\sigma$
for all $X,Y \in {\cal A}$.

\medskip

\begin{theorem} \cite[Theorem 7.1]{T:qRacah} \label{thm:dagger} \samepage
Let $(A; \{E_i\}_{i=0}^d; A^*; \{E^*_i\}_{i=0}^d)$
denote a Leonard system in $\cal A$.
Then there exists a unique antiautomorphism $\dagger$ of $\cal A$ such that
$A^\dagger=A$ and $A^{*\dagger}=A^*$.
Moreover $X^{\dagger\dagger}=X$ for all $X \in {\cal A}$.
\end{theorem}

\medskip

\begin{definition}          \label{def:D}         \samepage
Let $(A; \{E_i\}_{i=0}^d; A^*; \{E^*_i\}_{i=0}^d)$
denote a Leonard system in $\cal A$.
We let $\cal D$ denote the subalgebra of $\cal A$ generated by $A$.
We let ${\cal D}^*$ denote the subalgebra of $\cal A$ generated by $A^*$.
\end{definition}

\medskip

\begin{lemma}    \cite[Lemma 6.3]{T:survey}    \label{lem:dagD}    \samepage
Let $(A; \{E_i\}_{i=0}^d; A^*; \{E^*_i\}_{i=0}^d)$
denote a Leonard system in $\cal A$ and let
$\dagger$ denote the corresponding antiautomorphism of $\cal A$  from 
Theorem \ref{thm:dagger}.
Then referring to Definition \ref{def:D},
$\dagger$ fixes everything in $\cal D$ and everything in ${\cal D}^*$.
In particular 
\begin{equation}          \label{eq:dag}
     E^\dagger_i = E_i, \qquad  E^{*\dagger}_i = E^*_i
        \qquad (0 \leq i \leq d).
\end{equation}
\end{lemma}

\section{A basis for $\cal X$}

In this section we prove Theorem \ref{thm:main}.
We start with a lemma.

\medskip

\begin{lemma}    \cite[Lemma 11.1]{T:qRacah}          \samepage
Let $(A; \{E_i\}_{i=0}^d; A^*; \{E^*_i\}_{i=0}^d)$
denote a Leonard system in $\cal A$ and 
let $V$ denote a simple $\cal A$-module. Then 
$E_iV=E_iE^*_0V$ and $E^*_iV=E^*_iE_0V$
for $0 \leq i \leq d$.
\end{lemma}

\medskip

\begin{corollary}             \label{cor:YEiEs0}          \samepage
Let $(A; \{E_i\}_{i=0}^d; A^*; \{E^*_i\}_{i=0}^d)$
denote a Leonard system in $\cal A$.
Then for $Y \in {\cal A}$ the following hold for $0 \leq i \leq d$.
\begin{itemize}
\item[(i)] $YE_i=0$ if and only if $YE_iE^*_0=0$.
\item[(ii)] $YE^*_i=0$ if and only if $YE^*_iE_0=0$.
\end{itemize}
\end{corollary}

\medskip

\begin{corollary}             \label{cor:Es0EiY}          \samepage
Let $(A; \{E_i\}_{i=0}^d; A^*; \{E^*_i\}_{i=0}^d)$
denote a Leonard system in $\cal A$.
Then for $Y \in {\cal A}$ the following hold for $0 \leq i \leq d$.
\begin{itemize}
\item[(i)] $E_iY=0$ if and only if $E^*_0E_iY=0$.
\item[(ii)] $E^*_iY=0$ if and only if $E_0E^*_iY=0$.
\end{itemize}
\end{corollary}

\begin{proof}
Apply $\dagger$ to the equations in Corollary \ref{cor:YEiEs0},
and use Lemma \ref{lem:dagD}.
\end{proof}

\medskip

\begin{definition}           \label{def:theta}           \samepage
Let $(A; \{E_i\}_{i=0}^d; A^*; \{E^*_i\}_{i=0}^d)$
denote a Leonard system in $\cal A$.
For $0 \leq i \leq d$ we let $\theta_i$ (resp. $\theta^*_i$)
denote the eigenvalue of $A$ (resp. $A^*$)
associated with $E_i$ (resp. $E^*_i$).
We note that the scalars $\theta_0,\theta_1,\ldots,\theta_d$
(resp. $\theta^*_0,\theta^*_1,\ldots,\theta^*_d$) are mutually
distinct and contained in $\mathbb{K}$.
\end{definition}

\medskip

\begin{proposition}              \label{prop:key}            \samepage
Let $(A; \{E_i\}_{i=0}^d; A^*; \{E^*_i\}_{i=0}^d)$
denote a Leonard system in $\cal A$ and
let $\cal X$ denote the subspace of $\cal A$ from Definition \ref{def:X}.
Then for $X \in {\cal X}$ such that $XE^*_0=0$ and $XAE^*_0=0$ we have $X=0$.
\end{proposition}

\begin{proof}
First assume $d=0$. Then $E^*_0=I$ and
the result follows. For the rest of this proof assume $d\geq 1$.
We assume $X \neq 0$ and get a contradiction.

In the equation $I=\sum_{i=0}^d E^*_i$ we multiply each term on the right
by $AE^*_0$ and simplify the result using (\ref{eq:Atrid}) to obtain
$AE^*_0 = E^*_0AE^*_0 + E^*_1AE^*_0$;
expanding $XAE^*_0=0$ using this and $XE^*_0=0$ we find $XE^*_1AE^*_0=0$.
Let $V$ denote a simple $\cal A$-module and observe $XE^*_1AE^*_0V=0$.
Note that $E^*_1V= E^*_1AE^*_0V$, since $E^*_1AE^*_0 V \subseteq E^*_1V$,
$\dim E^*_1V=1$, and $E^*_1AE^*_0V\not=0$ in view of (\ref{eq:Atrid}).
By the above comments $XE^*_1V=0$ so $XE^*_1=0$.
In the equation $I=\sum_{i=0}^d E^*_i$ we multiply each term on the left 
by $E^*_0X$ and simplify the result using (\ref{eq:trids}) to find
$E^*_0X= E^*_0XE^*_0+E^*_0XE^*_1$;
now $E^*_0X=0$ since each of $XE^*_0$, $XE^*_1$ is zero. 

Since $X \neq 0$ there exist integers $i,j$ $(0 \leq i,j \leq d)$
such that $E_iXE_j \neq 0$. 
Define
\[
   r = \min\,\{ \min\{i,j\} \,|\, 0 \leq i,j \leq d,\; E_iXE_j \neq 0 \}.
\]

First assume $r=d$, so that $E_dXE_d\not=0$ and each of $E_iXE_d$, $E_dXE_i$ 
is zero for $0 \leq i \leq d-1$.
In the equation $I=\sum_{i=0}^d E_i$ we multiply each term on the left by $E_dX$
and simplify to get $E_dX= E_dXE_d$.
By this and since $XE^*_0=0$ we find $E_dXE_dE^*_0=0$.
Now $E_dXE_d=0$ by Corollary \ref{cor:YEiEs0}(i), for a contradiction.

Next assume $r\leq d-1$.
Note that for $0 \leq i \leq r-1$ we have
$E_rXE_i=0$ and $E_iXE_r=0$.
We now show that each of $E_rXE_r$ and $E_rXE_{r+1}$ is zero.
In the equation $I=\sum_{i=0}^d E_i$ we multiply each term on the left
by $E_r X$. We simplify the result using (\ref{eq:trid}) and our above
comments to find
\begin{equation}        \label{eq:keystar}
    E_r X = E_rXE_r + E_rXE_{r+1}.
\end{equation}
In this equation we multiply each term
on the right by $E^*_0$ and use $XE^*_0=0$ to find
\begin{equation}        \label{eq:key1}
    E_rX E_rE^*_0 + E_rX E_{r+1}E^*_0=0.
\end{equation}
We multiply each term of (\ref{eq:keystar})
on the right by $A$ and use 
$E_iA=\theta_iE_i$ $(0 \leq i \leq d)$ to find
$E_rXA = \theta_r E_rXE_r + \theta_{r+1}E_rXE_{r+1}$.
In this equation we multiply each term on the right by $E^*_0$
and use $XAE^*_0=0$ to find
\begin{equation}        \label{eq:key2}
  \theta_r E_rX E_rE^*_0 + \theta_{r+1} E_rX E_{r+1}E^*_0=0.
\end{equation}
Solving the linear system (\ref{eq:key1}) and (\ref{eq:key2}), we find
$E_rX E_rE^*_0=0$ and $E_rX E_{r+1}E^*_0=0$.
By this and Corollary \ref{cor:YEiEs0}(i) we find
$E_rX E_r=0$ and $E_rX E_{r+1}=0$.
Next we show $E_{r+1}XE_r=0$.
We mentioned earlier that $E_iXE_r= 0$ for $0 \leq i \leq r-1$.
In the equation $I= \sum_{i=0}^d E_i$ we multiply
each term on the right by $XE_r$. We simplify the result
using (\ref{eq:trid}) and our above comments to find
$XE_r = E_{r+1}XE_r$. In this equation we multiply each term on the
left by $E^*_0$ and use $E^*_0X=0$ to find
$E^*_0E_{r+1}XE_r=0$, so $E_{r+1}XE_r=0$ in view of 
Corollary \ref{cor:Es0EiY}(i).
We have now shown that each of
$E_rXE_r$, $E_rXE_{r+1}$, $E_{r+1}XE_r$ is zero,
contracting the  definition of $r$.
We conclude $X=0$.
\end{proof}

\medskip

\begin{corollary}             \label{cor:dimX}           \samepage
Let $(A; \{E_i\}_{i=0}^d; A^*; \{E^*_i\}_{i=0}^d)$
denote a Leonard system in $\cal A$.
Then the space $\cal X$ from Definition \ref{def:X} has
dimension at most $5$.
\end{corollary}

\begin{proof}
We assume $d \geq 2$; otherwise $\dim {\cal A} \leq 4$ and
the result follows.
We define linear maps 
$\pi_0:{\cal X} \to {\cal X}E^*_0$ and
$\pi_1:{\cal X} \to {\cal X}AE^*_0$ by
\[
  \pi_0(X)=XE^*_0, \qquad \pi_1(X)=XAE^*_0  \quad\qquad\qquad (X \in {\cal X}).
\]
For $i=0,1$ let $K_i$ denote the kernel of $\pi_i$.
We compute the dimensions of $K_0$ and  $K_1$. 
First observe
\[
   \dim E^*_i {\cal A} E^*_j = 1 \qquad   (0 \leq i,j\leq d).   
\]
We have ${\cal X} E^*_0 = E^*_0 {\cal X} E^*_0 + E^*_1 {\cal X} E^*_0$
in view of (\ref{eq:trids}); therefore $\dim {\cal X} E^*_0 \leq 2$ so
\begin{equation}   \label{eq:dimaux1}
       \dim K_0 \geq  \dim {\cal X} - 2.
\end{equation}
Combining (\ref{eq:Atrid}) and (\ref{eq:trids}) we routinely obtain
\[      
  {\cal X} A E^*_0 \subseteq
       E^*_0 {\cal A} E^*_0 + E^*_1 {\cal A} E^*_0 + E^*_2 {\cal A} E^*_0;
\]
therefore $\dim {\cal X} A E^*_0 \leq 3$ so
\begin{equation}          \label{eq:dimaux2}
        \dim K_1 \geq \dim {\cal X} - 3.
\end{equation}
The intersection of $K_0$ and $K_1$ is zero
by Propostion \ref{prop:key}; therefore
\begin{equation}          \label{eq:dimaux3}
    \dim K_0 + \dim K_1 \leq \dim {\cal X}.
\end{equation}
Combining (\ref{eq:dimaux1})--(\ref{eq:dimaux3}) we find 
$\dim {\cal X} \leq 5$ as desired.
\end{proof}

\medskip

\begin{proofof}{Theorem \ref{thm:main}}
Comparing (\ref{eq:Astrid}), (\ref{eq:Atrid}) and 
(\ref{eq:trid}), (\ref{eq:trids}) we see that
each of the elements (\ref{eq:basisX}) is contained in $\cal X$.
We must show they actually span $\cal X$, and that
they are linearly independent provided $d\geq 2$.
First assume $d=0$. Then the assertion is obvious.
Next assume $d=1$. Then one routinely verifies that 
${\cal X}={\cal A}$ is spanned by the elements  (\ref{eq:basisX}).
Finally assume $d \geq 2$.
In view of Corollary \ref{cor:dimX}, it suffices to show that
the elements (\ref{eq:basisX}) are linearly independent.
Suppose
\begin{equation}      \label{eq:spanaux}
e I + f A + f^* A^* + g AA^* + g^* A^*A=0
\end{equation}
for some scalars $e,f,f^*,g,g^*$ in $\mathbb{K}$.
We show each of $e,f,f^*,g,g^*$ is zero.
For $1 \leq i \leq d$ we multiply each term in (\ref{eq:spanaux})
on the left by $E^*_{i-1}$ and the right by $E^*_i$ to obtain
\[
   (f + g\theta^*_i + g^* \theta^*_{i-1})E^*_{i-1}AE^*_i=0.
\]
By this and since $E^*_{i-1}AE^*_{i}$ is nonzero we find
\begin{equation}            \label{eq:spanaux2}
   f + g\theta^*_i + g^* \theta^*_{i-1}=0   \qquad (1 \leq i \leq d).
\end{equation}
For $1 \leq i \leq d$ we multiply each term in (\ref{eq:spanaux})
on the left by $E^*_{i}$ and the right by $E^*_{i-1}$ to obtain
\[
   (f + g\theta^*_{i-1} + g^* \theta^*_{i})E^*_{i}AE^*_{i-1}=0.
\]
By this and since $E^*_{i}AE^*_{i-1}$ is nonzero we find
\begin{equation}            \label{eq:spanaux3}
   f + g\theta^*_{i-1} + g^* \theta^*_{i}=0   \qquad (1 \leq i \leq d).
\end{equation}
Combining (\ref{eq:spanaux2}) at $i=1$ and (\ref{eq:spanaux3}) 
at $i=1,2$ we routinely find that each of $f$, $g$, $g^*$ is zero.
Interchanging the roles of $A$ and $A^*$ in the above argument
we find $f^*=0$.
Now (\ref{eq:spanaux}) becomes $eI=0$ so $e=0$.
We have now shown that each of $e$, $f$, $f^*$, $g$, $g^*$ is zero 
and the result follows.
\end{proofof}

\section{The linear maps $\T$ and $\T^*$}

In this section we discuss some linear
maps $\T: {\cal X} \to {\cal D}$ and
$\T^* : {\cal X} \to {\cal D}^*$ that we find attractive.
To motivate things we recall some results
by the second author and Vidunas.

\medskip

\begin{lemma}   \cite[Theorem 1.5]{TV}   \label{lem:AW}       \samepage
Let $(A; \{E_i\}_{i=0}^d; A^*; \{E^*_i\}_{i=0}^d)$
denote a Leonard system in $\cal A$.
Then there exists a sequence of scalars
$\beta,\gamma,\gamma^*,\varrho,\varrho^*,\omega,\eta,\eta^*$
taken from $\mathbb{K}$ such that both
\begin{equation}        \label{eq:AW}
A^2A^* - \beta AA^*A+A^*A^2 - \gamma(AA^*+A^*A)-\varrho A^*
   = \gamma^* A^2 + \omega A + \eta I,
\end{equation}
\begin{equation}        \label{eq:AWs}
 A^{*2}A - \beta A^*AA^* + AA^{*2} - \gamma^*(A^*A+AA^*) -\varrho^* A
  = \gamma A^{*2} + \omega A^* + \eta^* I.
\end{equation}
Moreover the sequence is uniquely determined by 
the Leonard system provided $d \geq 3$.
\end{lemma}

\medskip

\begin{note}
The equations (\ref{eq:AW}) and (\ref{eq:AWs})
first appeared in \cite{Z}; they are called the
{\em Askey-Wilson relations}. 
\end{note}

\medskip

We have a comment.

\medskip

\begin{lemma}    \cite[Theorem 4.5]{TV}             \samepage
Let $(A; \{E_i\}_{i=0}^d; A^*; \{E^*_i\}_{i=0}^d)$
denote a Leonard system in $\cal A$.
Then referring to Definition \ref{def:theta} and Lemma \ref{lem:AW}
we have
\begin{eqnarray}
 \beta+1 &=&
   \frac{\theta_{i-2}-\theta_{i+1}}{\theta_{i-1}-\theta_{i}}
 =  \frac{\theta^*_{i-2}-\theta^*_{i+1}}{\theta^*_{i-1}-\theta^*_{i}}
          \qquad (2 \leq i \leq d-1),   \label{eq:beta} \\
  \gamma &=& \theta_{i-1} - \beta \theta_{i} + \theta_{i+1}
             \qquad (1 \leq i \leq d-1),    \label{eq:gamma}  \\
  \gamma^* &=& \theta^*_{i-1} - \beta \theta^*_{i} + \theta^*_{i+1}
             \qquad (1 \leq i \leq d-1),    \label{eq:gammas}  \\
 \varrho &=& \theta_{i-1}^2 - \beta \theta_{i-1}\theta_{i}+\theta_{i}^2
            -\gamma(\theta_{i-1}+\theta_{i})   \qquad (1 \leq i \leq d), 
                                             \label{eq:rho}  \\
 \varrho^* &=& \theta^{*2}_{i-1} - \beta \theta^*_{i-1}\theta^*_{i}
             + \theta^{*2}_{i} - \gamma^*(\theta^*_{i-1}+\theta^*_{i})   
                  \qquad (1 \leq i \leq d).          \label{eq:rhos} 
\end{eqnarray}
\end{lemma}

\medskip

\begin{theorem}          \label{thm:T}              \samepage
Let $(A; \{E_i\}_{i=0}^d; A^*; \{E^*_i\}_{i=0}^d)$
denote a Leonard system in $\cal A$.
Let the spaces $\cal X$ and $\cal D$ be as in Definitions \ref{def:X} and
\ref{def:D}, respectively.
Then there exists a $\mathbb{K}$-linear map
$\T : {\cal X} \to {\cal D}$ that satisfies
\begin{equation}        \label{eq:defT}
   \T(X)=A^2X-\beta AXA+XA^2-\gamma(AX+XA)-\varrho X
\end{equation}
for all $X \in {\cal X}$.
Moreover
\begin{eqnarray}
\T(I) &=& (2-\beta)A^2 - 2\gamma A-\varrho I,   \label{eq:TI}\\
\T(A) &=& (2-\beta)A^3 - 2 \gamma A^2 - \varrho A, \label{eq:TA} \\
\T(A^*) &=& \gamma^* A^2 + \omega A + \eta I,    \label{eq:TAs}  \\    
\T(AA^*) &=& \gamma^* A^3 + \omega A^2 + \eta A, \label{eq:TAAs} \\
\T(A^*A) &=& \gamma^* A^3 + \omega A^2 + \eta A. \label{eq:TAsA}
\end{eqnarray}
\end{theorem}

\begin{proof}
Certainly there exists a $\mathbb{K}$-linear map 
$\T:{\cal X} \to {\cal A}$ that satisfies (\ref{eq:defT}).
Using (\ref{eq:AW}) we find $\T$ satisfies
(\ref{eq:TI})--(\ref{eq:TAsA}).
Combining (\ref{eq:TI})--(\ref{eq:TAsA}) and Theorem \ref{thm:main}
we find $\T(X) \in {\cal D}$ for all $X \in {\cal X}$,
and the result follows.
\end{proof}

\medskip

Interchanging the roles of $A$ and $A^*$ in Theorem \ref{thm:T} 
we obtain:

\medskip

\begin{theorem}          \label{thm:Ts}              \samepage
Let $(A; \{E_i\}_{i=0}^d; A^*; \{E^*_i\}_{i=0}^d)$
denote a Leonard system in $\cal A$.
Let the spaces $\cal X$ and ${\cal D}^*$ be as in Definitions \ref{def:X} and
\ref{def:D}, respectively.
Then there exists a $\mathbb{K}$-linear map
$\T^* : {\cal X} \to {\cal D}^*$ that satisfies
\begin{equation}         \label{eq:defTs}
  \T^*(X)= A^{*2}X-\beta A^*XA^*+XA^{*2}-\gamma^*(A^*X+XA^*)-\varrho^* X
\end{equation}
for all $X \in {\cal X}$.
Moreover
\begin{eqnarray*}
\T^*(I) &=& (2-\beta)A^{*2} - 2\gamma^* A^* -\varrho^* I,  \\
\T^*(A^*) &=& (2-\beta)A^{*3} - 2 \gamma^* A^{*2} - \varrho^* A^*,
                                                \label{eq:TsAs} \\
\T^*(A) &=& \gamma A^{*2} + \omega A^* + \eta^* I, \label{eq:TsA} \\
\T^*(A^*A) &=& \gamma A^{*3} + \omega A^{*2} + \eta^* A, \label{eq:TsAsA} \\
\T^*(AA^*) &=& \gamma A^{*3} + \omega A^{*2} + \eta^* A.  \label{eq:TsAAs}
\end{eqnarray*}
\end{theorem}

\medskip

We have a comment concerning the image and kernel of $\T$.

\medskip

\begin{lemma}               \label{lem:kerT}     \samepage
Referring to Theorem \ref{thm:T} the following (i)--(iii) hold.
\begin{itemize}
\item[(i)] $\text{\rm Span}\{AA^*-A^*A\} \subseteq \text{\rm Ker}(\T)$.
\item[(ii)] $\text{\rm Im}(\T) \subseteq \text{\rm Span}\{I,A,A^2,A^3\}$.
\item[(iii)] Assume $d \geq 3$. Then equality holds in (i) if and only if 
equality holds in (ii).
\end{itemize}
\end{lemma}

\begin{proof}
(i), (ii): Immediate from Theorem \ref{thm:T}.

(iii): Use Theorem \ref{thm:main} and elementary linear algebra.
\end{proof}

\medskip

Interchanging the roles of $A$ and $A^*$ in Lemma \ref{lem:kerT} we obtain:

\medskip

\begin{lemma}               \label{lem:kerTs}          \samepage
Referring to Theorem \ref{thm:Ts} the following (i)--(iii) hold.
\begin{itemize}
\item[(i)] $\text{\rm Span}\{AA^*-A^*A\} \subseteq \text{\rm Ker}(\T^*)$.
\item[(ii)] $\text{\rm Im}(\T^*) \subseteq \text{\rm Span}\{I,A^*,A^{*2},A^{*3}\}$.
\item[(iii)] Assume $d \geq 3$. Then equality holds in (i) if and only 
if equality holds in (ii).
\end{itemize}
\end{lemma}

\medskip

Referring to Lemmas \ref{lem:kerT} and \ref{lem:kerTs} it appears that
we have equality in (i) and (ii) for most Leonard systems
but not all. Below we give an example where
equality is not attained.

\medskip

\begin{definition}        \label{def:bip}                 \samepage
Let $(A; \{E_i\}_{i=0}^d; A^*; \{E^*_i\}_{i=0}^d)$
denote a Leonard system in $\cal A$.
We say this Leonard system is {\em bipartite} (resp. {\em dual bipartite})
whenever  $E^*_iAE^*_i=0$ (resp. $E_iA^*E_i=0$)
for $0 \leq i \leq d$.
\end{definition}

\medskip

\begin{lemma}           \label{lem:bip}              \samepage
Let $\Phi=(A; \{E_i\}_{i=0}^d; A^*; \{E^*_i\}_{i=0}^d)$
denote a Leonard system in $\cal A$.
Then referring to Theorems \ref{thm:T}, \ref{thm:Ts} and Definition \ref{def:bip}
the following (i), (ii) hold provided $d \geq 3$.
\begin{itemize}
\item[(i)] 
Assume $\Phi$ is bipartite. Then
\begin{eqnarray*}
  \text{\rm Ker}(\T^*) &=& \text{\rm Span}\{A,AA^*,A^*A\},   \\
  \text{\rm Im}(\T^*) &=& \text{\rm Span}\{B^*, A^*B^* \},
\end{eqnarray*}
where $B^*=(2-\beta)A^{*2}-2\gamma^* A^* - \varrho^* I$.
\item[(ii)]
Assume $\Phi$ is dual bipartite. Then
\begin{eqnarray*}
  \text{\rm Ker}(\T) &=& \text{\rm Span}\{A^*,A^*A,AA^* \},   \\
  \text{\rm Im}(\T) &=& \text{\rm Span}\{ B, AB \},
\end{eqnarray*}
where $B=(2-\beta)A^2-2\gamma A- \varrho I$.
\end{itemize}
\end{lemma}

\begin{proof}
(ii):
By \cite{NT:balanced} and \cite[Theorem 5.3]{TV} 
each of $\gamma^*$, $\omega$, $\eta$ is zero. 
By this and Theorem \ref{thm:T} we have
$\text{\rm Ker}(\T) \supseteq \text{\rm Span}\{A^*, A^*A,AA^*\}$
and $\text{\rm Im}(\T) = \text{\rm Span}\{B, AB\}$.
To show 
$\text{\rm Ker}(\T) = \text{\rm Span}\{A^*, A^*A,AA^*\}$ 
it suffices to show that $B$ and $AB$ are linearly independent.
Suppose $B$ and $AB$ are linearly dependent.
Then $B=0$ since the elements $I,A,A^2,A^3$ are linearly independent.
Since $d\geq 3$ there exists an integer $i$ such that $1 \leq i \leq d-1$.
Multiplying each term in the equation 
$B=(2-\beta)A^2-2\gamma A - \varrho I$
by $E_i$ and simplifying we find $E_i$ times
\begin{equation}              \label{eq:bipaux1}
    (2-\beta)\theta_i^2 - 2 \gamma \theta_i - \varrho 
\end{equation}
is zero. Of course $E_i$ is not zero so (\ref{eq:bipaux1}) is zero.
Using (\ref{eq:gamma}) and (\ref{eq:rho}) we routinely find
(\ref{eq:bipaux1}) is equal to 
\[
    (\theta_i - \theta_{i-1})(\theta_i - \theta_{i+1})
\]
and is therefore nonzero.
This is a contradiction and the result follows.
\end{proof}

\medskip
\noindent
{\bf Open Problem:}
Referring to Lemmas \ref{lem:kerT} and \ref{lem:kerTs},
precisely determine the set of Leonard systems for which
equality holds in (i) and (ii).

\bigskip

\bibliographystyle{plain}

\bigskip\bigskip\noindent
Kazumasa Nomura\\
College of Liberal Arts and Sciences\\
Tokyo Medical and Dental University\\
Kohnodai, Ichikawa, 272-0827 Japan\\
email: nomura.las@tmd.ac.jp

\bigskip\noindent
Paul Terwilliger\\
Department of Mathematics\\
University of Wisconsin\\
480 Lincoln Drive\\ 
Madison, Wisconsin, 53706 USA\\
email: terwilli@math.wisc.edu

\bigskip\noindent
{\bf Keywords.}
Leonard pair, tridiagonal pair, $q$-Racah polynomial, orthogonal polynomial.

\noindent
{\bf 2000 Mathematics Subject Classification}.
05E35, 05E30, 33C45, 33D45.

\end{document}